\documentclass[11pt,oneside]{article}
\usepackage{color}
\usepackage{psfrag}
\usepackage{epsfig}
\usepackage{graphicx}
\usepackage{amsmath}
\usepackage{amssymb}
\usepackage{amsthm}
\usepackage{cite}
\usepackage[latin1]{inputenc}

\newtheorem{theorem}{Theorem}[section]

\newtheorem{corollary}[theorem]{Corollary}

\begin{document}

\title{Linear codes from Denniston maximal arcs}
\date{}
\author{Daniele Bartoli, Massimo Giulietti,   Maria Montanucci}

\maketitle

\begin{abstract}
In this paper we construct functional codes from Denniston maximal arcs. 
 For $q=2^{4n+2}$ we
obtain linear codes with parameters $[(\sqrt{q}-1)(q+1),5,d]_q$  where $\lim_{q \to +\infty} d=(\sqrt{q}-1)q-3\sqrt{q}$.
We also find  for $q=16,32$ a number of linear codes
which appear to have larger minimum distance with respect to the known codes  with same length and dimension.
\end{abstract}


\section{Introduction}

For $q$ a prime power, let $AG(N,q)$  denote the affine space over the finite field with $q$ elements $\mathbb F_q$. For a point set $\mathcal X\subseteq AG(N,q)$ and  a linear subspace  $\mathcal{V}\subset \mathbb{F}_q[X_1,\ldots, X_N]$,  the functional code $\mathcal{C}_{\mathcal{V}}(\mathcal{X}){\subset}\mathbb{F}_q^n$ is defined as
$$\mathcal{C}_{\mathcal{V}}(\mathcal{X}) := \{(f(P_1),\ldots,f(P_n)) \ : \ f \in \mathcal{V}\},$$
where  $\mathcal X=\{P_1,\ldots, P_n\}$. Clearly,
the code $\mathcal C_{\mathcal V}( \mathcal X)$ can be seen as  $\eta(\mathcal{V})$, where $\eta$ is the linear map
$$\eta\ : \ \mathcal{V} \to \mathbb{F}_q^n$$
with $\eta(f)=(f(P_1),\ldots,f(P_n))$ for any $f \in \mathcal{V}$. 
An analogous definition is given for $\mathcal X$ a subset of $PG(N,q)$, the $N$-dimensional projective space over $\mathbb F_q$; here, $\mathcal{V}$ consists of 
homogeneous  polynomials in $X_0,X_1,\ldots,X_N$.  
The length of  $\mathcal{C}_{\mathcal{V}}(\mathcal{X})$ is $n$, the dimension is $\dim_{\mathbb F_q} \eta(\mathcal{V})$, and the minimum distance coincides with 
$n-\max_{f \in \mathcal{V}, \, \eta(f)\neq 0} \#(V_f \cap \mathcal{X})$, where $V_f$ denotes the set of zeros of $f$.

The case where $\mathcal{X}$ consists of the set of $\mathbb F_q$-rational points of a quadric or a hermitian variety and $\mathcal{V}$ is a vector space of polynomials of a given degree has been thoroughly investigated; see for instance \cite{BDFS2014,FABE:2010A,FABE:2010B,AH:2010,BS2014,ELX2011}. In particular, functional codes from the Hermitian curve, the so-called Hermitain codes, have performances which are sometimes comparable with those of BCH codes; see e.g. {\cite[Ch. 4]{CJ2009} and \cite{JTZT2009, JC2006}.}

In this paper we construct functional codes arising from particular subsets $\mathcal{X} \in AG(2,q)$ (or $PG(2,q)$) called maximal $n$-arcs. A maximal $n$-arc $\mathcal X$  is  a set of $nq+n-q$ points  such that any line of the plane contains either $0$ or $n$ points of $\mathcal X$.  
The integer $n$ is called the degree of $\mathcal{X}$. In \cite{Barlotti55} Barlotti   {proved} that a necessary condition for the existence of a proper maximal arc is $n\mid q$, whereas it was shown in \cite{BBM1997} that  no nontrivial maximal arcs exist for $q$ odd.
Maximal arcs have interesting connection with linear codes; see for instance \cite{Delsarte1972,CK1986}. In particular, the functional code {$\mathcal C_{\mathcal V}(\mathcal X)$} where $\mathcal V$ is the set of linear forms and $\mathcal X$ is a maximal arc of degree $n$ in $PG(2,q)$ is an optimal $2$-weight code. 
Maximal arcs are also related to  partial geometries \cite{Mathon} and Steiner 2-designs \cite{DCVM1995,Thas1974}.

The classification of maximal arcs of degree $2$ in $PG(2,q)$,  also called hyperovals, is a long-standing and fascinating problem in Finite Geometry. The simplest example of hyperoval is given by  a conic plus its nucleus; for other infinite families see e.g. {\cite[Section 8.4]{HirschBook}}. {As to the higher degree case, in 1969 Denniston \cite{Denniston1969} gave a construction of maximal arcs of degree $n$ in Desarguesian projective planes of even order $q$, for all $n$ dividing $q$. Such maximal arcs consist of the union of some conics and their common nucleus.} Other constructions of maximal arcs are given in \cite{Thas1974Bis,Thas1980,Mathon,HM2004,Hamilton2002}; see also the references therein.

The aim of the paper is twofold. On the one hand, for small $q$'s we give explicit constructions of functional linear codes from Denniston maximal arcs having better parameters {than those listed in the}  database Mint \cite{Mint}; see Section \ref{DennistonCodes}. In particular, we obtain codes with parameters
$[119,10,94]_{16}$, $[120,5,103]_{16}$, $[119,5,103]_{16}$, $[51,5,43]_{16}$, $[99,5,88]_{32}$.

On the other hand, for $q=2^{4n+2}$ we show that some {functional codes from Denniston}  arcs reach the parameters of the functional Hermitian codes of the same dimension. We in fact  
obtain linear codes with parameters $[(\sqrt{q}-1)(q+1),5,d]_q$  where $\lim_{q \to +\infty} d=(\sqrt{q}-1)q-3\sqrt{q}$. 
This achievement relies on an interesting geometrical property of Denniston arcs. A Denniston maximal arc
can be seen as the set of $\mathbb F_q$-rational points of a (reducible) algebraic curve {of  degree $2t$ with equation $L(f(x,y))=0$, with $L$ an additive polynomial of degree $t$ and $f$ an irreducible quadratic form. For $q$ a square and $L(T)=T^{\sqrt{q}}-T$ we show that such curve intersects} a generic conic through the common nucleus of the conics in at most (roughly) $2\sqrt{q}$ points. This significantly improves the bound  $4\sqrt q$ provided by B\'ezout Theorem and shows that
 $\sqrt{q}(q-2)$ is a rough lower bound on the minimum distance of the functional code where $\mathcal V$ is the $5$-dimensional linear  space of polynomials generated by $X,Y,XY,X^2,Y^2$; see Theorem \ref{Th:Componenti} and Corollary \ref{Co:Componenti}.


\section{Codes from Denniston arcs of degree $q/2$ and $q/4$}\label{DennistonCodes}
We recall the construction of Denniston maximal arcs.
Let $H$ be an additive subgroup of $\mathbb{F}_q$. Consider an irreducible quadratic form $f(x,y)= ax^2+bxy+cy^2$  over $\mathbb{F}_q$.
Denniston maximal arcs{, viewed as subsets of $AG(2,q)$, can} described as 

$$\Omega =\{(x,y) \in AG(2,q) \ | \ f(x,y) \in H\};$$
see \cite{Denniston1969}. It is easily seen that $|\Omega|=n= (q+1)(|H|-1)+1$. {In this paper we will consider functional codes from both $\Omega$ and $\Omega^*=\Omega \setminus \{(0,0)\}$.} 

\subsection{$|H|=q/2$}
Consider now a subgroup $H$ of index $2$. Note that $(0,0)$ belongs to $\Omega$. Consider the $\mathbb{F}_q$-vector space $\mathcal{V}$ of dimension $5$ given by 
$$\{ax^2+bxy+cy^2+dx+ey \ : \ a,b,c,d,e \in \mathbb{F}_q\}.$$
Let $\eta : \mathcal{V} \to \mathbb{F}_q^{n}$ be defined by 
$$\eta(f) = (f(P_1),\ldots,f(P_n)),$$
where $\Omega=\{P_1,\ldots,P_n\}$. The image $\eta(\mathcal{V})$ is an $[n,5,n-\alpha]_q$ code, where $\alpha=\max_{f \in \mathcal{V}} \# \{V_f\cap \Omega\}$. In this case, if $f$ splits into two linear factors then it vanishes at  $q$ points of $\Omega$ at most. On the other hand, if $f$ defines an absolutely irreducible conic, $V_f$ contains at most $q+1$ points. So $\eta(\mathcal{V})$ is an $[(q^2-q)/2,5,(q^2-3q-2)/2]_q$ code.
Let us consider now $\Omega^*=\{P_1,\ldots,P_{n-1}\}$. By the same construction,  $\eta(\mathcal{V})$ is an $[(q^2-q)/2+1,5,(q^2-3q)/2+1]_q$ code, since a conic (irreducible or not) can contain at most $q$ points of the set $\Omega^*$. 

 Table \ref{Index2} shows the parameters of this family of codes. In particular, for $q=16$, the codes $[120,5,103]_{16}$ and $[119,5,103]_{16}$ improve the corresponding entries in the database Mint \cite{Mint}.

Now we deal with the case where $\mathcal{V}$ is the vector space of all the polynomials of degree at most three. By the Hasse-Weil bound {(see \cite[Theorem 5.2.3]{Sti})} an irreducible plane cubic curve has at most $16+1+2\sqrt{16}=25$ points in $PG(2,16)$. On the other hand, if the cubic splits in three lines or one line and a conic then it is easily seen that such a curve shares at most $24$ or $25$ with  $\Omega$ or $\Omega^*$. This ensures the existence of a $[120,10,95]_{16}$ and a $[119,10,94]_{16}$ which are better than the $[120,10,92]_{16}$ and the $[119,10,91]_{16}$ codes in \cite{Mint}.

\begin{table}
\caption{Intersection with conics, $|H| =q/2$}\label{Index2}
\begin{center}
\begin{tabular}{llll}
\hline
$q$ &$n$&$\eta(\mathcal{V})$&Best known\\
\hline
$16$&$120$&$[120,5,103]_{16}$&$[120,5,102]_{16}$\\
$16$&$119$&$[119,5,103]_{16}$&$[119,5,101]_{16}$\\
\hline

\end{tabular}
\end{center}
\end{table}

\subsection{$|H|=q/4$}
{We now consider the case where $\mathcal V$ is as in the previous subsection but $H$ is a subgroup of index $4$. As a result of a computer search, we found some functional codes on $\Omega^*$ over $\mathbb{F}_{16}$ and $\mathbb{F}_{32}$ with better parameters than those listed in \cite{Mint}. 
In Table 2  the weight distribution of these codes is described. Here $\eta$  is a primitive element of $\mathbb{F}_{16}$ satisfying $\eta^4 + \eta + 1=0$ and $\omega$ is a primitive element of $\mathbb{F}_{32}$ satisfying $\omega^5+\omega^2+1=0$. The group $H$ and the polynomial $f$ are specified in the second and in the third  column. The fourth column contains the weight distribution of $\eta(\mathcal{V})$. We remark that for $q=32$ several other different subgroups $H$ give rise to codes having the same weight distribution. An interesting open problem is to determine whether all these codes are isomorphic.}


\begin{table}
\caption{Codes with $|H|=4$}\label{Index4}
\tabcolsep = 1 mm
\begin{tabular}{|l|l|l|l|l|l|l|}
\hline
$q$& $H$&$f(x,y)$& Weight distribution& Parameters& Best known in \cite{Mint}\\
\hline
\hline

$16$& $\eta \mathbb{F}_{4}$&$x^2+ \eta^{10}x y+ \eta^8 y^2$&
\begin{tabular}{l}
$0^1 43^{459}44^{4272}45^{2992}$\\
$46^{5232}47^{12750}48^{18736}$\\
$49^{14280}50^{8496} 51^{2415}$\\
\end{tabular}
&\begin{tabular}{l}
$[51,5,43]_{16}$\\
\end{tabular}
&
\begin{tabular}{l}
$[51,5,42]_{16}$\\
\end{tabular}

\\

\hline
$32$& $\{0,\omega^9,\omega^{13},\omega^{19}\}$&$x^2+ x y+ y^2$&
\begin{tabular}{l}
$0^1 88^{66}89^{660}90^{1848}$\\
$91^{15774}92^{23628}93^{53592}$\\
$94^{110352}95^{197604}96^{251394}$\\
$97^{237732}98^{136488}99^{52206}$\\
\end{tabular}
&
\begin{tabular}{l}
$[99,5,88]_{32}$\\
\end{tabular}
&
\begin{tabular}{l}
$[99,5,87]_{32}$\\
\end{tabular}
\\

\hline
\end{tabular}

\end{table}

\section{Codes from Denniston arcs of degree $\sqrt{q}$}\label{Sec:Curva}

In this section we consider the case where $q=2^{4n+2}$ and the subgroup $H$ is  the field $\mathbb{F}_{\sqrt{q}}$. 
As $\sqrt q$ is not a power of $4$,  $\mathbb F_q=\mathbb F_{\sqrt q}(\xi)$, with $\xi$ a root of $T^2+T+1=0$.
Note that $\{1,\xi\}$ is a basis of $\mathbb{F}_{q}$ over $\mathbb{F}_{\sqrt{q}}$, and
the polynomial $T^2+T+\xi\in \mathbb{F}_q[T]$ is irreducible over $\mathbb{F}_q$. Then  
$$\Omega^* =\bigcup _{z \in \mathbb{F}_{\sqrt{q}}^*} \{(x,y) \in AG(2,q) \ | \ x^2+xy+\xi y^2 = z\}.$$

Consider a generic conic {$\mathcal{D}$} through $(0,0)$. Then the {affine} points $(\overline x,\overline y)$ in  {$\mathcal D$} can be parametrized as follows: 
\begin{equation}\label{parametrizzazione}
\overline x =\overline x(m)= -\frac{E +D m}{A+ Bm +C m ^2} \qquad \overline y=\overline y(m)= -m\frac{E+D m}{A+B m +C m ^2},
\end{equation}
where $A,B,C,D,E \in \mathbb{F}_{q}$ and $(A,B,C)\neq (0,0,0)$. 
{Note that values $m$ for which $A+ Bm +C m ^2=0$ correspond to ideal points  of the conic $\mathcal{D}$ and therefore we can suppose $A+ Bm +C m ^2\neq 0$.}

In order to determine the minimum distance of the functional code $\mathcal C_{\mathcal V}( \Omega^*)$, with 
$\mathcal V$ the linear space generated by $X,Y,X^2,XY,Y^2$,
we need to count the possible intersections between {$\mathcal{D}$}  and $\Omega^*$. This is equivalent to determine the number of pairs $(m,z) \in \mathbb{F}_q\times \mathbb{F}_{\sqrt{q}}$ such that 
\begin{equation}\label{CondizioneIniziale}
\left(-\frac{E +D m}{A+ Bm +C m ^2}\right)^2(1+m+\xi m^2)= z.
\end{equation}

Write $m=m_1+\xi m_2$, $A=a_1+\xi a_2$, $B=b_1+\xi b_2$, $C=c_1+\xi c_2$, $D=d_1+\xi d_2$, $E=e_1+\xi e_2$, with $m_i,a_i,b_i,c_i,d_i,e_i\in \mathbb{F}_{\sqrt{q}}$ for $i=1,2$. Then Equation \eqref{CondizioneIniziale} reads
\begin{equation}\label{Condizione}
\left\{
\begin{array}{l}
f(m_1,m_2,z)=0\\
g(m_1,m_2,z)=0\\

\end{array}
\right.,
\end{equation}
where 
$$
\begin{array}{l}
f(m_1,m_2,z)=(a_1 + a_2 + b_1 m_1 + b_1 m_2 + b_2 m_1 + c_1 m_1^2 + c_2 m_1^2 + 
        c_2 m_2^2)^2 z + d_1^2  m_1^3\\  
        \qquad + d_1^2  m_1^2 + d_1^2  m_1  m_2^2 + d_1^2  m_2^4 + d_1^2  m_2^3 + d_1^2  m_2^2 
        + d_2^2  m_1^4 + d_2^2  m_1^3 
        + d_2^2  m_1^2  m_2 \\ 
        \qquad + d_2^2  m_1^2+  d_2^2  m_2^3 + e_1^2  m_1 + e_1^2  m_2^2 + e_1^2 + e_2^2  m_1^2 + e_2^2  m_1 
        + e_2^2  m_2^2 + e_2^2  m_2 + e_2^2,\\
        \\
g(m_1,m_2,z)=(a_2 + b_1 m_2 + b_2 m_1 + b_2 m_2 + c_1 m_2^2 + c_2 m_1^2)^2z  + d_1^2  m_1^4  \\
        \qquad + d_1^2  m_1^2  m_2+ d_1^2  m_1  m_2^2 + 
        d_1^2  m_2^4 + d_1^2  m_2^2 + d_2^2  m_1^3 + d_2^2  m_1^2 + d_2^2  m_1  m_2^2 
           \\
        \qquad + d_2^2  m_2^4+ d_2^2  m_2^3+ d_2^2  m_2^2+ e_1^2  m_1^2 + e_1^2  m_2 + 
        e_2^2  m_1 + e_2^2  m_2^2 + e_2^2.\\
\end{array}
$$


{The coefficients of $z$ in $f$ and $g$ cannot be both the zero polynomial in $m_1$ and $m_2$ otherwise $(A,B,C)=(0,0,0)$.}
{First of all we note that the bijection $(m_1,m_2,z) \mapsto (m_1,m_2,z+1)$ does not change the number of triples $(M_1,M_2,Z) \in \mathbb{F}_{\sqrt{q}}^3$ satisfying  \eqref{Condizione}. Let $f^{\prime}(m_1,m_2,z)=f(m_1,m_2,z+1)$ and $g^{\prime}(m_1,m_2,z)=f(m_1,m_2,z+1)$. Also, both  $f^{\prime}$ and $g^{\prime}$ are linear in $z$. We denote by $h(m_1,m_2)$ the polynomial obtained by eliminating $z$. 
In the proof of Theorem \ref{Th:Componenti} we will show that the number of triples $(M_1,M_2,Z) \in \mathbb{F}_{\sqrt{q}}^3$ satisfying \eqref{Condizione} is at most the number of pairs $(M_1,M_2) \in \mathbb{F}_{\sqrt{q}}^2$ satisfying $h(M_1,M_2)=0$.  Note that by a straightforward computation $h$ can be written as  $\sum_{0\leq i,j\leq 8} {\alpha_{i,j}}m_1^im_2^j$ with}

{
\begin{equation}\label{Pol_h}
{\normalsize
\begin{array}{lll}
\alpha_{8,0}=\alpha_{0,8}=c_2^2& 
\alpha_{4,4}=c_2^2&
\alpha_{7,0}=c_1^2\\
\\
\alpha_{6,2}=\alpha_{2,6}=c_2^2&
\alpha_{5,2}=\alpha_{3,4}=c_1^2 + c_2^2&
\alpha_{2,5}=\alpha_{4,3}= c_1^2\\
\\
\alpha_{0,7}=c_1^2 + c_2^2&
\alpha_{6,0}=\begin{array}{l}b_2^2 + c_1^2 e_1^2 + c_1^2\\ + c_2^2 e_1^2 + c_2^2 e_2^2\\ \end{array}&
\alpha_{4,2}=\alpha_{2,4}=\begin{array}{l}b_1^2 + b_2^2 + c_1^2 e_2^2\\ + c_1^2 + c_2^2 e_1^2 + c_2^2\\ \end{array}\\
\\
\alpha_{0,6}=\begin{array}{l}b_1^2 + c_1^2 e_1^2 + c_1^2 e_2^2\\ + c_2^2 e_2^2 + c_2^2\\ \end{array}&
\alpha_{5,0}=b_1^2 + c_1^2 e_2^2 + c_2^2 e_1^2&
\alpha_{4,1}=  b_2^2 + c_1^2 e_1^2 + c_2^2 e_1^2 + c_2^2 e_2^2\\
\\
\alpha_{3,2}=b_1^2&
\alpha_{2,3}= b_2^2&
\alpha_{1,4}=b_1^2 + c_1^2 e_1^2 + c_1^2 e_2^2 + c_2^2 e_2^2\\
\\
\alpha_{0,5}=b_2^2 + c_1^2 e_2^2 + c_2^2 e_1^2&
\alpha_{4,0}=\begin{array}{l}
a_2^2 + b_1^2 e_1^2 + b_1^2 + b_2^2 e_1^2\\ + b_2^2 e_2^2 + c_1^2 e_2^2 + 
c_2^2 e_1^2\\\end{array}&
\alpha_{2,2}=b_1^2 e_1^2 + b_1^2 + b_2^2 e_1^2 + b_2^2 e_2^2\\
\\
\alpha_{0,4}=\begin{array}{l}a_1^2 + a_2^2 + b_1^2 e_1^2 + b_1^2\\ + b_2^2 e_1^2 + b_2^2 e_2^2 + c_1^2 e_1^2 \\
+ c_1^2 e_2^2 + c_2^2 e_2^2\\\end{array}&
\alpha_{3,0}=a_1^2 + b_1^2 e_2^2 + b_2^2 e_1^2&
\alpha_{2,1}=\begin{array}{l}a_2^2 + b_1^2 e_1^2 + b_2^2 e_1^2\\ + b_2^2 e_2^2\\ \end{array}\\
\\
\alpha_{1,2}=\begin{array}{l}a_1^2 + a_2^2 + b_1^2 e_1^2\\ + b_2^2 e_1^2 + b_2^2 e_2^2\\ \end{array}&
\alpha_{0,3}=a_1^2 + b_1^2 e_1^2 + b_1^2 e_2^2 + b_2^2 e_2^2&
\alpha_{2,0}=\begin{array}{l}a_1^2 e_1^2 + a_1^2 + a_2^2 e_1^2\\ + a_2^2 e_2^2 + b_1^2 e_2^2 + b_2^2 e_1^2\\\end{array}\\
\\
\alpha_{0,2}=\begin{array}{l}a_1^2 e_2^2 + a_1^2 + a_2^2 e_1^2 + a_2^2\\ + b_1^2 e_1^2 + b_2^2 e_1^2 + 
b_2^2 e_2^2\\\end{array}&
\alpha_{1,0}=a_1^2 e_2^2 + a_2^2 e_1^2&
\alpha_{0,1}=a_1^2 e_1^2 + a_2^2 e_1^2 + a_2^2 e_2^2\\
\\
\alpha_{0,0}=a_1^2 e_2^2 + a_2^2 e_1^2.\\
\end{array}}
\end{equation}
}

The polynomial $h(m_1,m_2)$ defines a plane curve $\mathcal{X}$ of order {at most} $8$. {First we show that $\mathcal{X}$ has at most two absolutely irreducible components defined over $\mathbb{F}_{\sqrt{q}}$. By the Hasse-Weil Theorem, this will give us an upper bound on the number of solutions of $h(m_1,m_2)=0$, and hence on the number or triples $(M_1,M_2,Z) \in \mathbb{F}_{\sqrt{q}}^3$ satisfying \eqref{Condizione}.} 

Note that $\mathcal{X}$ contains the points 
$$P_1=(\eta,\eta^2), \quad P_2=(\eta^2,\eta^4), \quad P_3=(\eta^4,\eta^8), \quad P_4=(\eta^8,\eta)$$
 with $\mathbb{F}_{16}^*=\langle \eta\rangle$, {where  $\eta^4 + \eta + 1=0$.} We distinguish a number of cases.
 
\begin{enumerate}
\item $D\neq 0$. We can suppose $(d_1,d_2)=(0,1)$. The homogeneous part of degree $8$ of $h(m_1,m_2)$ is $c_2^2(m_1^2+m_1m_2+m_2^2)^4$.
\begin{enumerate}
\item $c_2\neq 0$. Since the ideal points of $\mathcal{X}$ are not $\mathbb{F}_{\sqrt{q}}$-rational, there are no $\mathbb{F}_{\sqrt{q}}$-rational lines contained in $\mathcal{X}$. It is easy to see that no $\mathbb{F}_{\sqrt{q}}$-rational cubic component can be contained in $\mathcal{X}$. Also, conic or quartic components of $\mathcal{X}$ are $\mathbb{F}_{\sqrt{q}}$-rational if and only if their homogeneous part of highest degree is $(m_1^2+m_1m_2+m_2^2)$ or $(m_1^2+m_1m_2+m_2^2)^2$ respectively. If $\mathcal{X}$ has more that two absolutely irreducible components defined over $\mathbb{F}_{\sqrt{q}}$ then it must split in: 2 lines and 3 conics, 4 conics, 2 conics and 1 quartic. 

In the first case at least two points among $P_1,P_2,P_3,P_4$ are not contained in the two lines, which must be of type $m_1 + \xi m_2+\alpha=0$ and $m_1 + \xi^2m_2+\beta=0$. So at least one conic must contain a point $P_i$. By direct checking this implies that the conic is not defined over $\mathbb{F}_{\sqrt{q}}$.

In the second case, a  conic of equation $m_1^2+m_1m_2+m_2^2+Am_1+Bm_2+C=0$ contains all the points $P_i$'s if and only if 
$$
\left\{
\begin{array}{l} 
\eta^4 A + \eta^8 B + C + \eta^3=0\\
\eta A + \eta^2 B + C + \eta^{12}=0\\
\eta^8 A + \eta B + C+ \eta^6=0\\
\eta^2 A + \eta^4 B + C + \eta^9=0\\
\end{array}
\right..
$$
The previous system has no solution. Therefore there exist at least two conics which contain the points $P_i$'s. On the other hand, by direct checking and recalling that $\{1,\eta,\eta^5,\eta^6\}$ is a basis of $\mathbb{F}_{q^2}$ over $\mathbb{F}_{\sqrt{q}}$,  such conics are not defined over $\mathbb{F}_{\sqrt{q}}$. This implies that at most two conics can be $\mathbb{F}_{\sqrt{q}}$-rational.

In the third case, arguing as above, we conclude that the points $P_i$'s must be contained in the quartic component. The quartic $\mathbb{F}_{\sqrt{q}}$-rational component $\mathcal{Q}$ must be defined by
$$(m_1^2+m_1  m_2+m_2^2)^2+A_1  m_1^3+A_2  m_1^2  m_2+A_3  m_1  m_2^2+A_4  m_2^3$$
$$\qquad \qquad+B_1  m_1^2+B_2  m_1  m_2+B_3  m_2^2+C_1  m_1+C_2  m_2+C_3=0,$$
where $A_i,B_i,C_i \in \mathbb{F}_{\sqrt{q}}$. The condition $P_i \in \mathcal{Q}$ yields 
$$A_1 = B_1 + B_2 + C_1 + C_2 + C_3, \quad A_2 = B_3 + C_3, \quad A_3 = C_1 + C_3 + 1, \quad A_4 = B_1 + C_1 + C_2 + C_3 + 1.$$

Also, if $\mathcal{Q}$ and two conics {$\mathcal{D}_1$} and {$\mathcal{D}_2$} of equation $m_1^2+m_1  m_2+m_2^2+\alpha_i m_1+\beta_i m_2+\gamma_i=0$ are all components of $\mathcal{X}$, then 
$$B_3 + C_1 + 1=0, \qquad B_2 + B_3 + C_3 + 1=0.$$

Now, by direct checking the point $Q = (\eta^5 e_1 + \eta^{10} e_2 + \eta^3,\eta^{10} e_1 + e_2  + \eta^{13})$ belongs to $\mathcal{X}$. If $Q\in {\mathcal{D}_i}$, then the conic splits into two non-$\mathbb{F}_{\sqrt{q}}$-rational lines, so $Q$ must be contained in $\mathcal{Q}$. This yields $B_1=C_1=0$, $B_3=C_2=C_3=1$ and in this case $\mathcal{Q}$ splits in four lines not defined over $\mathbb{F}_{\sqrt{q}}$.

\item $c_2=0$ and $c_1\neq 0$. In this case $\mathcal{X}$ has degree $7$. The homogeneous part of degree 7 is given by $c_1^2(m_1^2+m_1m_2+m_2^2)^3(m_1+m_2)$. An $\mathbb{F}_{\sqrt{q}}$-rational line contained in $\mathcal{X}$ can only have equation $m_1+m_2+A_1=0$. By direct checking this implies 
$$
\left\{
\begin{array}{l}
b_1 + b_2 + c_1e_2 + c_1A_1=0\\
c_1(b_1 + b_2 + c_1 e_2 + \eta^2 c_1)(b_1 + b_2 + c_1 e_2 + \eta c_1)=0\\
\end{array}
\right..
$$
This is not possible since $c_1\neq 0$ by assumption. Therefore the only cases in which $\mathcal{X}$ has more than two absolutely irreducible components defined over $\mathbb{F}_{\sqrt{q}}$ are: 1 line and 3 conics, 2 conics and 1 cubic. 

In the first case, the line $\ell$ must be of type $m_1+m_2+A_1=0$, otherwise not all the conics are $\mathbb{F}_{\sqrt{q}}$-rational. Arguing as above, we immediately notice that not all the points $P_i$'s can be contained in $\ell$, and this forces at least one conic to be non-$\mathbb{F}_{\sqrt{q}}$-rational.

In the second case, all the points $P_i$'s and the point $Q$ must be contained in the cubic which has equation 
$$(m_1^2+m_1m_2+m_2^2)(m_1+m_2)+A_1m_1^2+A_2 m_1m_2+A_3m_2^2+B_1m_1+B_2 m_2+B_3=0.$$
This yields
$$A_1 + B_2 + 1=0,\qquad  A_3=B_1=B_3,\qquad  A_2=0, \qquad e_2^2+e_2+1=0,$$
impossibile since $e_2 \in \mathbb{F}_{\sqrt{q}}$. 

\item $c_1=c_2=0$ and $b_2\neq 0$. In this case $\mathcal{X}$ has degree 6. The homogeneous part of degree 6 is given by $(b_1 m_2 + b_2 m_1)^2(m_1^2+m_1m_2+m_2^2)^2$. A linear  $\mathbb{F}_{\sqrt{q}}$-rational component of $\mathcal{X}$ should have equation $b_1 m_2 + b_2 m_1+A_1=0$. This implies $b_1=b_2$, $a_1 = b_2 e_2$, and $(b_2 + \eta A_1)(b_2 + \eta^2 A_1)=0$, impossible. The unique case in which we have more than two absolutely irreducible components defined over $\mathbb{F}_{\sqrt{q}}$ is given by two conics of type $m_1^2+m_1m_2+m_2^2+Am_1+Bm_2+C=0 $ and one conic of type $(b_1 m_2 + b_2 m_1)^2+Am_1+Bm_2+C=0$. Suppose that all the  $P_i$'s and $Q$ belong to the last conic. By direct checking such a conic has equation $b_2^2 m_1^2 + b_1^2 m_2^2 + b_1^2 m_1    + b_2^2 m_2 +b_1^2 =0$ and both 
$$b_1^2  e_2^2 + b_1^2  e_2 + b_1^2 + b_2^2  e_1^2 + b_2^2  e_1 + b_2^2  e_2 + b_2^2=0$$
and 
$$b_1^2  e_1^2 + b_1^2  e_1 + b_1^2  e_2 + b_1^2 + b_2^2  e_1^2 + b_2^2  e_1 + b_2^2  e_2^2=0.$$
This gives $b_2=0$, impossible. So at least one conic of equation $(b_1 m_2 + b_2 m_1)^2+Am_1+Bm_2+C=0$ contains one point among $\{P_1,P_2,P_3,P_4,Q\}$ and, as above, it is not $\mathbb{F}_{\sqrt{q}}$-rational. 

\item $c_1=c_2=b_2=0$ and $b_1\neq 0$. In this case $\mathcal{X}$ has degree 6 and the homogeneous part of degree 6 is given by $b_1^2m_2^2(m_1^2+m_1m_2+m_2^2)^2$. An $\mathbb{F}_{\sqrt{q}}$-rational line contained in $\mathcal{X}$ should have equation $m_2=A_1$. By direct checking this is impossible since $b_1\neq 0$. As in the previous case the unique decomposition of $\mathcal{X}$ that we need to check is given by three conics. Since one of them should have equation $m_2^2+Am_1+Bm_2+C=0$, if $P_i$'s and $Q$ belong to it, then $A=C=1$ and $B=0$ and $e_2^2 + e_2 + 1=0$, impossible. So at least one conic of equation $(b_1 m_2 + b_2 m_1)^2+Am_1+Bm_2+C=0$ contains one point among $\{P_1,P_2,P_3,P_4,Q\}$ and, as above, it is not $\mathbb{F}_{\sqrt{q}}$-rational.

\item $c_1=c_2=b_1=b_2=0$. In this case $\mathcal{X}$ has degree 4 and the homogeneous part of degree 4 is given by $(\alpha_1m_2 + \alpha_2 m_1 + \alpha_2 m_2)^4$, where $\alpha_i^2=a_i$, $i=1,2$. A linear component of $\mathcal{X}$ should be $\alpha_1m_2 + \alpha_2 m_1 + \alpha_2 m_2 +A=0$. It is easily seen that such a line cannot be a component of $\mathcal{X}$ and therefore the number of $\mathbb{F}_{\sqrt{q}}$-rational components of $\mathcal{X}$ is at most two.

\end{enumerate}

\item $D= 0$. We can suppose $(e_1,e_2)=(1,0)$, since otherwise the conic {$\mathcal{D}$} splits into two lines. In this case $\mathcal{X}$ has degree 6 and the homogeneous part of degree 6 is given by $(c_1 m_1 + c_1 m_2 + c_2 m_1)^2(m_1^2+m_1m_2+m_2^2)^2$. 

\begin{enumerate}
\item $(c_1,c_2)\neq (0,0)$. It is easily seen that there exists no linear $\mathbb{F}_{\sqrt{q}}$-rational component in $\mathcal{X}$. 

If $e_2=0$, $e_1\neq 0$, $b_1=b_2=a_1=1$ then $\mathcal{X}$ splits in 
$$e_1^2(m_1 + m_2 + \eta^5)(m_1 + m_2 + \eta^{10})(a_2^2 + m_1^4 + m_1^3 + m_1^2 m_2^2 + m_1 m_2 + m_2^4 + m_2^3 + m_2^2 + m_2)$$
and it contains at most two conic components defined over $\mathbb{F}_{\sqrt{q}}$. Otherwise, $\mathcal{X}$ does not contain linear components. The only possibility is three conics. If one of the points $P_i$'s belongs to a conic of equation $m_1^2+m_1 m_2+m_2^2+\alpha_i m_1+\beta_im_2+\gamma_i=0$ then, as above, it is not $\mathbb{F}_{\sqrt{q}}$-rational. So all the points $P_i$'s must belong to the conic of equation $(c_1 m_1 + c_1 m_2 + c_2 m_1)^2+Am_1+Bm_2+C=0$ which is hence of type
$$c_1^2  m_1^2 + c_1^2  m_1 + c_1^2  m_2^2 + c_1^2  m_2 + c_1^2 + c_2^2  m_1^2 + c_2^2  m_2=0.$$
The other two conics must be of equation  $m_1^2+m_1 m_2+m_2^2+\alpha_i m_1+\beta_im_2+\gamma_i=0$, $i=1,2$. Easy computations show that in this case $c_1=c_2=0$, impossible.

\item $(c_1,c_2)= (0,0)$. If $b_1\neq b_2$ then $\mathcal{X}$ has degree 4 and the homogeneous part of degree 4 is given by $(b_1+b_2)^2(m_1^2+m_1m_2+m_2^2)^2$. No $\mathbb{F}_{\sqrt{q}}$-rational linear components are contained in $\mathcal{X}$ and therefore it contains at most two absolutely irreducible components defined over $\mathbb{F}_{\sqrt{q}}$. 

If $b_1=b_2$ then $\mathcal{X}$ has degree 3 and the homogeneous component of degree 3 is given by $(b_1)^2(m_1+m_2)(m_1^2+m_1m_2+m_2^2)$.  Clearly, at most one $\mathbb{F}_{\sqrt{q}}$-rational linear components is contained in $\mathcal{X}$ and therefore it contains at most two absolutely irreducible component defined over $\mathbb{F}_{\sqrt{q}}$. 

\end{enumerate}
\end{enumerate}

{We are now in a position to prove the main result of this section.} 

\begin{theorem}\label{Th:Componenti}
{The size of the intersection between a conic {$\mathcal{D}$} through the origin and the set $\Omega^*$ is at most $2q+2+20\sqrt{q}$. }
\end{theorem}
\proof
{As shown above, the number of absolutely irreducible $\mathbb{F}_{\sqrt{q}}$-rational components of the curve $\mathcal{X}$ is two. By the Hasse-Weil Theorem (see \cite[Theorem 5.2.3]{Sti}) it is easily seen that the maximum number of affine $\mathbb{F}_{\sqrt{q}}$-points $(M_1,M_2)$ of $\mathcal{X}$ is 
$$2q+20\sqrt{q}+2.$$}

{Let $r_1(m_1,m_2)$ and $r_2(m_1,m_2)$ be given by  
$$
\begin{array}{lll}
r_1(m_1,m_2)&=&a_1 + a_2 + b_1 m_1 + b_1 m_2 + b_2 m_1 + c_1 m_1^2 + c_2 m_1^2 +c_2 m_2^2,\\
r_2(m_1,m_2)&=&a_2 + b_1 m_2 + b_2 m_1 + b_2 m_2 + c_1 m_2^2 + c_2 m_1^2.\\
\end{array}
$$
We do not have to deal with pairs $(M_1,M_2)\in \mathbb{F}_{\sqrt{q}}^2$ such that $r_1(M_1,M_2)=r_2(M_1,M_2)=0$, since otherwise $M=M_1+\xi M_2$ is such that $A+BM+CM^2=0$ and this would give an ideal point of $\mathcal{D}$, see the parametrization of $\mathcal{D}$ in \eqref{parametrizzazione}. }

{So, for a pair $(M_1,M_2)\in \mathbb{F}_{\sqrt{q}}^2$ such that $h(M_1,M_2)=0$, we have that  one between $r_1(M_1,M_2)$  and $r_2(M_1,M_2)$ is different from $0$. This yields the existence of a unique $Z\in \mathbb{F}_{\sqrt{q}}$ such that the triple $(M_1,M_2,Z)$ satisfies the  conditions in \eqref{Condizione}. 
}
\endproof

\begin{corollary}\label{Co:Componenti} For $q=2^{4n+2}$,
there exist linear codes with parameters $[(\sqrt{q}-1)(q+1),5,d]_{q}$,  with $$d\geq (\sqrt{q}-1)(q+1)-(2q+2+20\sqrt{q})= (\sqrt{q}-3)q-19\sqrt{q}-3.$$
\end{corollary}

{\proof Using the notation of this section, the code $\mathcal{C}_{\mathcal{V}}(\Omega^*)$, where 
$$\Omega^* =\bigcup _{z \in \mathbb{F}_{\sqrt{q}}^*} \{(x,y) \in AG(2,q) \ | \ x^2+xy+\xi y^2 = z\}$$
and $\mathcal{V}=\langle X,Y,X^2,XY,Y^2\rangle$, has length equal to $|\Omega^*|=(\sqrt{q}-1)(q+1)$. Its dimension is $\dim(\mathcal{V})=5$ and by Theorem \ref{Th:Componenti} its minimum distance is at least $(\sqrt{q}-1)(q+1)-(2q+2+20\sqrt{q})=(\sqrt{q}-3)q-19\sqrt{q}-3$.
\endproof
}

\section{Concluding remarks and open questions}

In this paper we constructed functional codes arising from Denniston maximal arcs. In some cases such codes have better parameters than the already known ones. Some computations have been done using the software MAGMA \cite{Magma}. Due to computational reasons, we could  perform the search only for small $q$. In these cases, we observed that the weight distribution of the codes associated with different subgroups were the same. We do not know if different subgroups of the same size can give rise to inequivalent codes or not. 

For $q=2^{4n+1}$ we find linear codes with dimension $5$, length with the same order of magnitude of $q^3$ and 
Singleton defect bounded by $2q$. Functional codes with similar parameters can be obtained from the
 Hermitian curve of $PG(2,q)$. We were not able to compare the weight distributions of the two codes.

\section*{Acknowledgment}
The  authors were supported in part by Ministry for Education, University and Research of Italy (MIUR) (Project PRIN 2012 ``Geometrie di Galois e strutture di incidenza") and by the Italian National Group for Algebraic and Geometric Structures and their Applications (GNSAGA - INdAM).

\begin{flushleft}
Daniele Bartoli\\
Department of Mathematics and Computer Science,\\
University of Perugia,\\
Italy\\
e-mail: {\sf daniele.bartoli@unipg.it}
\end{flushleft}

\begin{flushleft}
Massimo Giulietti\\
Department of Mathematics and Computer Science,\\
University of Perugia,\\
Italy\\
e-mail: {\sf massimo.giulietti@unipg.it}
\end{flushleft}

\begin{flushleft}
Maria Montanucci\\
Dipartimento di Matematica Informatica ed Economia,\\
Universit\`a degli Studi della Basilicata,\\
Italy\\
e-mail: {\sf maria.montanucci@unibas.it}
\end{flushleft}

\end{document}